\documentclass[a4paper,12pt,openany]{article}
\usepackage[english]{babel}
\usepackage{amsmath,amsbsy,amssymb,latexsym,epsfig,graphicx}
\begin{document}
\title{An inverse of the modular invariant}
\author{S. Adlaj\footnote{\textit{e-mail: SemjonAdlaj@gmail.com}}}
\date{May 30, 2011}
\parskip 6pt
\maketitle
\begin{abstract}
During the last few years of his life, Ramanujan had adamantly tried to invert the modular invariant. Subsequent efforts failed until May 30, 2011 when an explicit closed formula for an inverse was presented at the CCRAS (Moscow, Russia). This very formula, along with some special values of the modular invariant, is given in this paper. 
\end{abstract}

\noindent
In a previous paper \cite{A}, a justification for defining an \textit{essential elliptic function} was made. Yet, enabling an inversion of the modular invariant is, perhaps, even more convincing. We shall not elaborate upon describing previous attempts for inverting the modular invariant aside from mentioning two typical references \cite{B,C}. The first reference provides a glimpse upon Ramanujan latest efforts, whereas the appendix of the second concludes with a well-known expression for a point $\tau$ in the fundamental domain as a ratio of hypergeometric functions, thereby linking $\tau$ with an intermediate variable $\lambda$. Formula (3.3), in the same paper, yields the modular invariant $j$ as a (well-known) fractional transformation of $\lambda$, of degree 6. We point out this transformation so as to suggest that verifying a formula for an inverse of the modular invariant is as straightforward as verifying a root of a given hexic.

\noindent 
An inversion of the modular invariant is afforded via successively composing the functions
$$k_0(x) = \frac{i \hspace{.04cm} G \left( \sqrt{1 - x^2} \ \right)}{G(x)}, \ k_1(x) = \frac{\sqrt{x + 4} - \sqrt{x}}{2}, \ k_2(x) = \frac{3}{2} \left( \frac{x}{k_3(x)} + k_3(x) \right) - 1,$$
where
$$k_3(x) = \sqrt[3]{\sqrt{x^2 - x^3} - x}$$
and $G(x)$ is the arithmetic-geometric mean of $1$ and $x$.
%and $G(\cdot)$ is the Gauss hypergeometric function satisfying the differential equation
%$$4 x \left( x - 1 \right) G''(x) + 4 \left( 2 x - 1 \right) G'(x) + G(x) \equiv 0.$$
In other words, the function
$$k = k_0 \circ k_1 \circ k_2$$
is an inverse of the modular invariant, which (we need not point out) is not single-valued.

\subsection*{An ascending sequence of special values of the modular invariant on the boundary of the fundamental domain and the \mbox{imaginary axis}}
$$j \left( 1/2 + 2 \sqrt{- 1} \ \right) = \frac{\left( \displaystyle \frac{4}{3} \left( 1 - 2 \left( 33 + 24 \sqrt{2} - 4 \sqrt{140 + 99 \sqrt{2}} \ \right)^2 \right)^2 - 1 \right)^3}{\left( 1 - 2 \left( 33 + 24 \sqrt{2} - 4 \sqrt{140 + 99 \sqrt{2}} \ \right)^2 \right)^2 - 1} <$$
$$< j \left( 1/2 + \sqrt{- 3} \ \right) = \frac{\left( \displaystyle 12 \left( 555 + 16 \left( 20 \sqrt{3} - \left( 8 / \sqrt{3} + 5 \right) \sqrt{ 26 + 15 \sqrt{3}} \ \right) \right)^2 - 1 \right)^3}{9 \left( 555 + 16 \left( 20 \sqrt{3} - \left( 8 / \sqrt{3} + 5 \right) \sqrt{ 26 + 15 \sqrt{3}} \ \right) \right)^2 - 1} <$$
$$< j \left( 1/2 + \sqrt{- 2} \ \right) = \frac{\left( \displaystyle \frac{4}{3} \left( 1 - 2 \left( 5 + 4 \sqrt{2} - 2 \sqrt{2 \left(7 + 5 \sqrt{2} \ \right)} \ \right)^2 \right)^2 - 1 \right)^3}{\left( 1 - 2 \left( 5 + 4 \sqrt{2} - 2 \sqrt{2 \left( 7 + 5 \sqrt{2} \ \right) } \ \right)^2 \right)^2 - 1} <$$
$$< j \left( 1/2 + \sqrt{- 1} \ \right) = \left( 181 - 19 \left( \frac{3}{\sqrt{2}} \right)^3 \right)^3 < j \left( \frac{1 + \sqrt{- 3}}{2} \ \right) = 0 <$$
$$< j \left( \frac{1 + 2\sqrt{- 2}}{3} \ \right) = \left( \frac{5 \left( 19 - 13\sqrt{2} \ \right)}{6} \right)^3 < j \left( \frac{1 + 4 \sqrt{- 3}}{7} \ \right) =$$
$$ \frac{\left( \displaystyle 12 \left( 555 - 16 \left( 20 \sqrt{3} - \left( 8 / \sqrt{3} - 5 \right) \sqrt{ 26 - 15 \sqrt{3}} \ \right) \right)^2 - 1 \right)^3}{9 \left( 555 - 16 \left( 20 \sqrt{3} - \left( 8 / \sqrt{3} - 5 \right) \sqrt{ 26 - 15 \sqrt{3}} \ \right) \right)^2 - 1} < j \left( \sqrt{- 1} \ \right) = 1 <$$
<$$< j \left( \frac{2}{\sqrt{-3}} \ \right) = \frac{375 \left( 35010 - 20213 \sqrt{3} \ \right)}{16} < j \left( \sqrt{- 2} \ \right) = \left( \frac{5}{3} \right)^3 < j \left( \sqrt{- 3} \ \right) = \frac{125}{4} <$$
$$< j \left( 2 \sqrt{- 1} \ \right) = \left( \frac{11}{2} \right)^3 < j \left( \frac{4}{\sqrt{- 3}} \ \right) =$$
$$= \frac{\left( \displaystyle 12 \left( 555 - 16 \left( 20 \sqrt{3} + \left( 8 / \sqrt{3} - 5 \right) \sqrt{ 26 - 15 \sqrt{3}} \ \right) \right)^2 - 1 \right)^3}{9 \left( 555 - 16 \left( 20 \sqrt{3} + \left( 8 / \sqrt{3} - 5 \right) \sqrt{ 26 - 15 \sqrt{3}} \ \right) \right)^2 - 1} < j \left( 2\sqrt{- 2} \ \right) =$$
$$= \left( \frac{5 \left( 19 + 13\sqrt{2} \ \right)}{6} \right)^3 < j \left( 2\sqrt{- 3} \ \right) =  \frac{375 \left( 35010 + 20213 \sqrt{3} \ \right)}{16} < j \left( 4 \sqrt{- 1} \ \right) =$$
$$= \left( 181 + 19 \left( \frac{3}{\sqrt{2}} \right)^3 \right)^3 < j \left( 4 \sqrt{- 2} \right) = \frac{\left( \displaystyle \frac{4}{3} \left( 1 - 2 \left( 5 + 4 \sqrt{2} + 2 \sqrt{ 2 \left( 7 + 5 \sqrt{2} \ \right) } \right)^2 \right)^2 - 1 \right)^3}{\left( 1 - 2 \left( 5 + 4 \sqrt{2} + 2 \sqrt{2 \left( 7 + 5 \sqrt{2} \ \right)} \right)^2 \right)^2 - 1} <$$
$$< j \left( 4 \sqrt{- 3} \right) = \frac{\left( \displaystyle 12 \left( 555 + 16 \left( 20 \sqrt{3} + \left( 8 / \sqrt{3} + 5 \right) \sqrt{ 26 + 15 \sqrt{3}} \ \right) \right)^2 - 1 \right)^3}{9 \left( 555 + 16 \left( 20 \sqrt{3} + \left( 8 / \sqrt{3} + 5 \right) \sqrt{ 26 + 15 \sqrt{3}} \ \right) \right)^2 - 1} <$$
$$< j \left( 8 \sqrt{- 1} \right) = \frac{\left( \displaystyle \frac{4}{3} \left( 1 - 2 \left( 33 + 24 \sqrt{2} + 4 \sqrt{140 + 99 \sqrt{2}} \ \right)^2 \right)^2 - 1 \right)^3}{\left( 1 - 2 \left( 33 + 24 \sqrt{2} + 4 \sqrt{140 + 99 \sqrt{2}} \ \right)^2 \right)^2 - 1}.$$

\begin{figure}[p]
\begin{picture}(500,520)(-2,0)
\put(210,616){\parbox{330pt}{$4 \sqrt{-1}$}}
\put(310,300){\parbox{50pt}{$2 \sqrt{-1}$ \scriptsize{$\displaystyle + \ \frac{1}{2}$}}}
\put(310,260){\parbox{44pt}{$\sqrt{-3}$ \scriptsize{$\displaystyle + \ \frac{1}{2}$}}}
\put(310,210){\parbox{44pt}{$\sqrt{-2}$ \scriptsize{$\displaystyle + \ \frac{1}{2}$}}}
\put(210,136){\parbox{330pt}{$\sqrt{-1}$}}
\put(225,-20){\parbox{36pt}{\normalsize{$0$}}}
\put(-10,-20){\parbox{36pt}{\scriptsize{$\displaystyle -\frac{3}{2}$}}}
\put(450,-20){\parbox{36pt}{\scriptsize{$\displaystyle \frac{3}{2}$}}}
\put(105,-20){\parbox{36pt}{\scriptsize{$\displaystyle -\frac{\sqrt{2}}{2}$}}}
\put(128,-20){\parbox{36pt}{\scriptsize{$\displaystyle -\frac{\sqrt{3}}{3}$}}}
\put(310,-20){\parbox{36pt}{\scriptsize{$\displaystyle \frac{\sqrt{3}}{3}$}}}
\put(328,-20){\parbox{36pt}{\scriptsize{$\displaystyle \frac{\sqrt{2}}{2}$}}}
\put(116,-60){The fundamental domain with some points,}
\put(48,-80){at which the value of the modualr invariant is calculated, being marked}
%\put(168,-140){SemjonAdlaj@gmail.com}
\includegraphics[width=160mm]{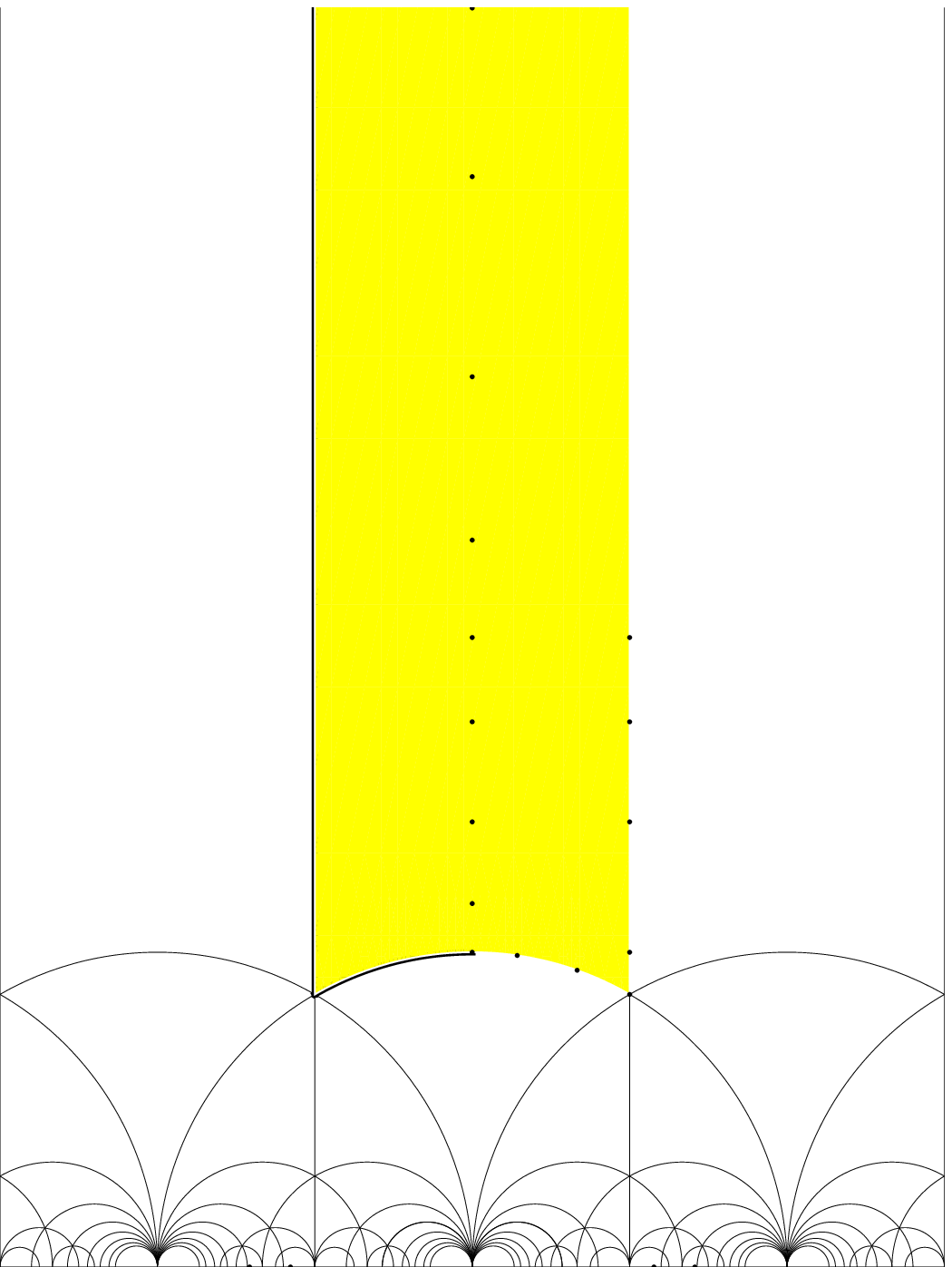}
\end{picture}
\end{figure}

\end{document}